\documentclass[reqno]{amsart}

\usepackage{amsmath}
\usepackage{amsfonts}

\newtheorem{thm}{Theorem}[section]
\newtheorem{lem}{Lemma}[section]
\newtheorem{conj}{Conjecture}[section]

\theoremstyle{definition}
\newtheorem{defin}{Definition}[section]

\def\fibo#1#2{\Biggl(\binom{#1}{#2}\Biggr)}

\begin{document}
\title{The Filbert Matrix}
\vspace{.5in}

\author[T. M. Richardson]{Thomas M. Richardson \\
Department of Mathematics and Statistics\\
Western Michigan University\\
Kalamazoo, MI  49008}

\email{ribby@alumni.umich.edu}

\begin{abstract}
A Filbert matrix
is a matrix whose $(i,j)$ entry is
$1/F_{i+j-1}$,
where $F_n$ is the $n^{\rm th}$ Fibonacci number.
The inverse of the $n\times n$ Filbert matrix
resembles the inverse of the $n\times n$ Hilbert
matrix, 
and we prove that it shares the property of
having integer entries.
We prove that the matrix formed by
replacing the Fibonacci numbers with
the Fibonacci polynomials has entries which
are integer polynomials.
We also prove that certain Hankel matrices
of reciprocals of binomial coefficients
have integer entries, and we conjecture
that the corresponding matrices based
on Fibonomial coefficients have integer
entries.
Our method is to give explicit formulae
for the inverses.
\end{abstract}

\maketitle

\section{Introduction}
The $n \times n$ Hilbert matrix is the $n \times n$
matrix whose
$(i,j)$-entry is $\frac{1}{i+j-1}$.
In \cite{TRICKS}, Man-Duen Choi explores
many fascinating properties of the Hilbert matrix,
including the fact the the $(i,j)$-entry of its
inverse is
\begin{equation}
\label{hilInvEnt}
\alpha_{ij}=
(-1)^{i+j}
(i+j-1)
\binom{n+i-1}{n-j}
\binom{n+j-1}{n-i}
\binom{i+j-2}{i-1}^2.
\end{equation}
Choi asks what sort of coincidence it is 
if the inverse of a matrix of reciprocals of
integers
has integer entries.
In this paper we show that
the inverses of the Hankel matrices based on the
reciprocals of the Fibonacci numbers,
the reciprocals of the binomial coefficients
$\binom{i+j}{2}$, and
the reciprocals of the binomial coeffiencts
$\binom{i+j+2}{3}$
all have integer entries.
We also find formulas for the entries of the
inverses of
these matrices and related matrices.

\begin{defin}
Let $\{a_k\}$ be an integer sequence 
with $a_k\ne 0$ for $k\ge 1$.
A {\em reciprocal Hankel matrix } 
based on $\{a_k\}$
is a matrix whose
$(i,j)$-entry is $1/a_{i+j-1}$.
We denote the $n\times n$ reciprocal Hankel
matrix based on $\{a_k\}$ by $R_n(a_k)$.
\end{defin}

The formula for the entries of the inverse of
$R_n(F_k)$ bears a striking resemblence to the
formula for the entries of the inverse of the
Hilbert matrix. 
Therefore, we call
a reciprocal Hankel matrix based on
the Fibonacci numbers a {\em Filbert matrix}.

\section{Filbert matrices}

We need the Fibonomial coefficents 
to describe the inverse of the Filbert
matrix.
See \cite{TAOCP} for more information on the
Fibonomial 
coefficients.

\begin{defin}
The {\em Fibonomial coefficients} are 
$$
\Biggl(\binom{n}{k}\Biggr)=\prod_{i=1}^k\frac{F_{n-i+1}}{F_i},
$$
where $n$ and $k$ are nonnegative integers.
\end{defin}

\begin{thm}
Let
$e(n,i,j)=n(i+j+1)+\binom{i}{2}+\binom{j}{2}+1$,
and let
$W(n)$ be the $n\times n$ matrix whose
$(i,j)$-entry is
$$
W_{ij}(n)=
(-1)^{e(n,i,j)}
F_{i+j-1}
\Biggl(\binom{n+i-1}{n-j}\Biggr)
\Biggl(\binom{n+j-1}{n-i}\Biggr)
\Biggl(\binom{i+j-2}{i-1}\Biggr)^2
.$$
Then the $n\times n$ matrix $W(n)$ is
the inverse of the Filbert matrix $R_n(F_k)$,
and $W(n)$ is an integer matrix.
\end{thm}

This theorem is a special case of Theorem
\ref{FilbertX},
which we prove below.
The formula for the entries of the inverse
closely corresponds to the formula
for the entries of the inverse of the $n\times n$
Hilbert matrix.
It results from Eq.\eqref{hilInvEnt} by 
changing all binomial coeffiecients to Fibonomial
coefficients
and changing the exponent of $-1$.
The pattern of the signs of entries the inverse of
$R_n(F_k)$ is
that they are constant on $2\times 2$ blocks, and
alternate
between blocks.

The Fibonacci polynomials $f_n(x)$ are defined by
$f_0(x)=0$, $f_1(x)=1$, and
$f_n(x)=xf_{n-1}(x)+f_{n-2}(x)$
for $n\ge 2$.
We also use $f_n$ to denote the Fibonacci
polynomial $f_n(x)$,
especially when we want to reduce the clutter in
some equations.
The $x$-Fibonomial coefficients are
the obvious generalization of 
the Fibonomial coefficients.

\begin{defin}
The {\em $x$-Fibonomial coefficients} are
$$
\fibo{n}{k}_x=\prod_{i=1}^k\frac{f_{n-i+1}(x)}{f_i(x)}
,$$
where $n$ and $k$ are nonnegative integers.
\end{defin}

To form the $(i,j)$-entry of the inverse of
$R_n(f_k(x))$,
replace each Fibonacci number and Fibonomial
coefficient
in $W_{ij}(n)$ with the corresponding Fibonacci
polynomial
and $x$-Fibonomial coefficient.

\begin{thm}
\label{FilbertX}
Let $V(n)$ be the $n\times n$ matrix whose
$(i,j)$-entry is
$$
V_{ij}(n)=
(-1)^{e(n,i,j)}
f_{i+j-1}
\Biggl(\binom{n+i-1}{n-j}\Biggr)_x
\Biggl(\binom{n+j-1}{n-i}\Biggr)_x
\Biggl(\binom{i+j-2}{i-1}\Biggr)_x^2
.$$
Then the $n\times n$ matrix $V(n)$ is
the inverse of the Filbert matrix $R_n(f_k(x))$,
and the entries of $V(n)$ are integer polynomials.
\end{thm}

The recurrence
$
\Bigl(\binom{n}{k}\Bigr)_x=
f_{k-1}(x)\Bigl(\binom{n-1}{k}\Bigr)_x+
f_{n-k+1}(x)\Bigl(\binom{n}{k}\Bigr)_x
$
shows that the Fibonomial coeffcients are
integer polynomials,
which implies that the entries of $V(n)$
are integer polynomials.

\section{Technology}
The proof of Theorem \ref{FilbertX}
and proofs of succeeding theorems
amount to proving various identities involving
sums of products of 
Fibonomial coefficients and binomial coefficients.
We supply computer proofs of these identities.
In some cases, the computer cannot do the
entire proof directly, and human intervention
is required to separate the proof into smaller
pieces that can be done by computer.

The programs and packages used to produce
the proofs for this paper include 
Maple V Release 5,
the Maple package {\tt EKHAD}
written by Doron Zeilberger,
and the Mathematica package {\tt MultiSum}
written by Kurt Wegschaider. 
{\tt EKHAD} is described in
\cite{AeqB}, and it is available at {\tt
www.math.temple.edu/$\sim$zeilberg}.
{\tt MultiSum} 
is available through the web site
{\tt www.risc.uni-linz.ac.at/software/},
and it 
is described in \cite{Weg}.
The particular functions that we use
from these packages are
{\tt zeil} from {\tt EKHAD}
and
{\tt FindRecurrence} from {\tt MultiSum}.

Both of these functions find a telescoped 
recurrence for a summand $F(n,k)$, where
$k$ is the summation variable.
The function {\tt zeil} uses Zeilberger's
algorithm to find a rational function $R(n,k)$
and a recurrence operator $P(n,N)$, 
where $N$ is the shift operator in $n$,
such that
\begin{equation}
\label{zeq}
P(n,N)(F(n,k))=R(n,k+1)F(n,k+1)-R(n,k)F(n,k).
\end{equation}
Let $f(n)$ be the unrestricted sum $\sum_kF(n,k)$. 
In many situations, 
Eq. \eqref{zeq} implies that $P(n,N)f(n)=0$,
making it easy to verify that $f(n)$ is constant.

The function {\tt FindRecurrence} gives similar
results with summands of the form
$F({\bf n},{\bf k})$ where ${\bf n}$ and ${\bf k}$
are vectors. 

Maple V Release 5 also includes an implementation
of Zeilberger's algorithm as the function
{\tt sumrecursion} of the package {\tt sumtools}.
However, {\tt sumrecursion}
only gives the recurrence operator $P(n,N)$,
and not the rational function $R(n,k)$,
which will be essential when we prove identities
involving a restricted sum. 

The sums involved in the
proof of  Theorem \ref{FilbertX}
are of products of Fibonomials, not
binomials, so these procedures
do not apply.
However, we obtained recurrences
for sums of products of Fibonomials
by modifying recurrences
found by these procedures
for the corresponding
sums of products of binomials.
 
\section{Proof of Theorem \ref{FilbertX}}

The $(i,m)$-entry of the product
$V(n)R_n(f_k(x))$
is
\begin{equation*}
p(n,i,m)=\sum_{j=1}^n P(n,i,m,j),
\end{equation*}
where
\begin{multline*}
\lefteqn{P(n,i,m,j)=}\\
(-1)^{e(n,i,j)}
\frac{f_{i+j-1}}{f_{j+m-1}}
\Biggl(\binom{n+i-1}{n-j}\Biggr)_x
\Biggl(\binom{n+j-1}{n-i}\Biggr)_x
\Biggl(\binom{i+j-2}{i-1}\Biggr)_x^2.
\end{multline*}
The summand satisfies the following 
recurrence relation that is
related to a recurrence produced by
{\tt FindRecurrence}
for an entry of the product of the Hilbert matrix
and its inverse.

\begin{lem}
\label{filRecLem}
The summand $P(n,i,m,j)$ satisfies the recurrence
relation
\begin{multline}
\label{filRecEqn}
{}-f_{n-i+1}f_{n+i-2}\Bigl(P(n,i-1,m,j)-P(n-1,i-1,m,j)\Bigr)+
\\
(-1)^{n+i}f_{i-1}^2\Bigl(P(n,i,m,j)-P(n-1,i,m,j)\Bigr)=0,
\end{multline}
and the sum $p(n,i,m)$ satisfies the recurrence
relation
\begin{multline}
\label{filSumEqn}
{}-f_{n-i+1}f_{n+i-2}\Bigl(p(n,i-1,m)-p(n-1,i-1,m)\Bigr)+
\\
(-1)^{n+i}f_{i-1}^2\Bigl(p(n,i,m)-p(n-1,i,m)\Bigr)=0.
\end{multline}
\end{lem}

\begin{proof}
Write each of the terms in Eq. \eqref{filRecEqn}
as a multiple of $P(n-1,i-1,m,j)$
to get the equation
\begin{multline}
-f_{n-i+1}f_{n+i-2}\bigl(P(n,i-1,m,j)-P(n-1,i-1,m,j)\bigr)
\\
(-1)^{n+i}f_{i-1}^2\bigl(P(n,i,m,j)-P(n-1,i,m,j)\bigr)
\\
=\frac{f_{n+i-2}}
{f_{n-i+1}f_{n-j}f_{i+j-1}}
M(n,i,j)
P(n-1,i-1,m,j)
\end{multline}
where
\begin{multline}
M(n,i,j)
=
(-1)^{i+j}
f_{n+i-1}f_{n+j-1}f_{i+j-2}
+
f_{n-i}f_{n-j}f_{i+j-2}
\\
 + (-1)^{i+j-1}
f_{n+i-2}f_{n+j-1}f_{i+j-1}
+
f_{n-i+1}f_{n-j}f_{i+j-1}.
\end{multline}
It suffices to show that $M(n,i,j)=0$.
But this follows from the standard Fibonacci
identities
$f_{n-i}f_{i+j-2}+f_{n-i+1}f_{i+j-1}=f_{n+j-1}$
and
$f_{n+i-2}f_{i+j-1}-f_{n+i-1}f_{i+j-2}=(-1)^{i+j-2}f_{n-j}$.

\end{proof} 

If we can establish 
$p(n,1,1)=1$,
$p(n,1,m)=0$ if $m\ne 1$, and
$p(n,n,n)=1$,
then Eq. \eqref{filSumEqn} shows that
$p(n,i,m)=1$ if $i=m$ and
$p(n,i,m)=0$ if $i\ne m$,
for $1\le i,m\le n$.

Case $p(n,1,m)$:
The summand $P(n,1,m,j)$
satisfies the recurrence
\begin{multline}
\label{pn1mEqn}
(-1)^{m+1}f_{n-1}f_{n+m-2}P(n,1,m-1,j)-f_nf_{n-m+1}P(n-1,1,m-1,j)
\\
+(-1)^{m}f_{n-1}f_{n+m-1}P(n,1,m,j)+f_nf_{n-m}P(n-1,i,m,j)=0,
\end{multline}
and this implies a similar recurrence for
$p(n,1,m)$.
The proof of Eq. \eqref{pn1mEqn}
is similar to the proof
of Lemma \ref{filRecLem}.
The initial values of this recurrence are 
$p(m,1,m)$ and $p(n,1,1)$.
The summand $P(m,1,m,j)$ satisfies the recurrence
$$
(-1)^mf_mf_{m-1}P(m,1,m,j)=G_1(m,j+1)-G_1(m,j)
$$
where
$G_1(m,j)=(-1)^{j-1}f_jf_{j-1}P(m,1,m,j)$.
Since the support of $G_1$ is $2\le j \le m$,
this equation implies that
$(-1)^mf_mf_{m-1}p(m,1,m)=0$.
Therefore, when $m>1$ we get $p(m,1,m)=0$.
Finally, the summand $P(n,1,1,j)$ satisfies 
$$
(-1)^nf_n^2P(n,1,1,j)=G_2(n,j+1)-G_2(n,j),
$$
where
$
G_2(n,j)=(-1)^{j-1}f_j^2P(n,1,1,j)
$.
In this case, the support of $G_2$ is $1\le j \le
n$,
so summing over $j$ from 1 to $n$ gives
$(-1)^nf_n^2p(n,1,1)=-G_2(n,1)=(-1)^nf_n^2$,
implying $p(n,1,1)=1$.

Case $p(n,n,n)$:
The summand $P(n,n,n,j)$
satisfies the recurrence
\begin{equation*}
P(n+1,n+1,n+1,j)-P(n,n,n,j)=
G_3(n,j+1)-G_3(n,j),
\end{equation*}
where
\begin{equation*}
G_3(n,j)=
(-1)^{e(n,n,j)}
\Biggl(
\frac{f_{3n+j-1}}{f_{n+j-1}}+2(-1)^n
\Biggr)
\Biggl(\binom{2n-1}{n-j+1}\Biggr)_x
\Biggl(\binom{n+j-2}{j-2}\Biggr)_x^2.
\end{equation*}
When we sum over $j$, the right hand side
telescopes to 0
and the left side is $p(n+1,n+1,n+1)-p(n,n,n)$.

This completes the proof of Theorem \ref{FilbertX}.

\section{Reciprocal Hankel matrices based on
binomial coefficients}

In this section we will prove that certain
reciprocal matrices
based on binomial coefficients have integer
entries.
We will give formulas for the entries of the
inverses
of these matrices.

Let $a_k=\binom{k+1}{2}$.

\begin{thm}
Let $A(n)$ be the $n\times n$ matrix whose
$(i,j)$-entry is
$$
A_{ij}(n)=\sum_{k=0}^{j-1}
(-1)^{i+k+1}
\binom{n+i}{n-k}\binom{n+k}{n-i}
\binom{i+k-1}{k}\binom{i+k}{k}
\frac{i}{2}.
$$
Then $A_{ij}(n)$ is an integer,
and $A(n)$ is
the inverse of the matrix $R_n(a_k)$.
\end{thm}

\begin{proof}
First we show that $A_{ij}(n)$ is an integer.
We use the well known fact that if $a$ is even
and $b$ is odd, then $\binom{a}{b}$ is even.
If $i$ is even, then obviously $A_{ij}(n)$ is an
integer,
so assume that $i$ is odd.
Now if $k$ is also odd, then $\binom{i+k}{k}$ is
even,
so we may assume that $k$ is even. Now one of
$\binom{n+i}{i+k}$ and $\binom{n+k}{i+k}$ is even.

Theorem \ref{Mother} below
shows that that $A(n)$ is
the inverse of the matrix $R_n(a_k)$.
\end{proof}

Let $b_k=b_k(r)$ be the binomial coefficient
$\binom{k+r-1}{r}$.
Suppose that $r$ a positive integer and $r\ge 3$.
Then the inverse of $R_n(b_k(r))$ does not always
have integer entries, but the values of $n$ for
which
the inverse does have integer entries seem to occur
periodically. Further, when the entries are not
integers,
the denominators are divisors of $r$.
The following conjecture is true for 
$n\le 20$, $r\le 10$, and $r$ an integer. 

\begin{conj}
Suppose that $r$ is a positive integer.
The inverse of the matrix $R_n(b_k(r))$
has integer entries if and only if
$n\equiv 0 \pmod q$ or $n\equiv 1 \pmod q$
for all prime powers $q$ that divide $r$.
\end{conj}

We do have an explicit formula for the entries of
the inverse.

\begin{thm}
\label{Mother}
Let $B(n,r)$ be the $n\times n$ matrix whose
$(i,j)$-entry is
\begin{multline*}
B_{ij}(n,r)=\\
\sum_{k=0}^{j-1}
(-1)^{i+k+1}
\binom{n+i+r-2}{i}\binom{n}{i}
\binom{n+k+r-2}{k}\binom{n}{k}
\frac{i^2\prod_{l=0}^{r-3}i+j+l}
{r\prod_{l=0}^{r-2}i+k+l}
.
\end{multline*}
Then $B(n,r)$ is the inverse of the matrix
$R_n(b_k)$.
\end{thm}

The theorem is valid if $r$ is an
indeterminate, not just if it is a positive
integer.
Also note that 
$B_{ij}(n,1)$ simplifies to $\alpha_{ij}$, 
the $(i,j)$-entry of the inverse of the Hilbert
matrix,
and $B_{ij}(n,2)$ is equal to $A_{ij}(n)$.

\begin{proof}
Let 
\begin{multline*}
H(n,i,m,j,k)=(-1)^{i+k+1}
\binom{n+i+r-2}{i}
\binom{n}{i}
\times
\\
\binom{n+k+r-2}{k}
\binom{n}{k}\frac{i^2\prod_{l=0}^{r-3}i+j+l}
{r\prod_{l=0}^{r-2}i+k+l}
\frac{1}{\binom{j+m+r-2}{r}},
\end{multline*}
so that
$h(n,i,m)=\sum_{j=1}^n\sum_{k=0}^{j-1}
H(n,i,m,j,k)$
is the $(i,m)$-entry of
$B(n,r)R_n(b_k)$.
Then $H$ satisfies the recurrence
\begin{equation}
\begin{split}
n^2(i-m+r-1)(n-i+r-1)(n+i+r-3)H(n-1,i-1,m-1,j,k)
&-\\
n^2(i-m-1)(n-i+r-1)(n+i+r-3)H(n-1,i-1,m,j,k)
&+\\
n^2(i-1)^2(i-m+1)H(n-1,i,m-1,j,k) &-\\
n^2(i-1)^2(i-m-r+1)H(n-1,i,m,j,k) &-\\
(n+r-2)^2(i-m+r-1)(n-i+1)(n+i-1)H(n,i-1,m-1,j,k)
&+\\
(n+r-2)^2(i-m-1)(n-i+1)(n+i-1)H(n,i-1,m,j,k)
&-\\
(n+r-2)^2(i-1)^2(i-m+1)H(n,i,m-1,j,k)
&+\\
(n+r-2)^2(i-1)^2(i-m-r+1)H(n,i,m,j,k) &=0.
\end{split}
\end{equation}
The preceeding recurrence was found by
{\tt FindRecurrence}.
The theorem will follow if we can establish the
correct values of
$h(n,1,m)$, $h(n,n,n)$, and $h(n,i,1)$.

Case $h(n,1,m)$:
Maple computes $h(n,1,1)=1$,
and it computes
\begin{equation*}
H_1(n,1,m,j)=\sum_{k=0}^{j-1}H(n,1,m,j,k)=
\frac
{(-1)^{j+1}j\binom{n+j+r-2}{j}\binom{n}{j}}
{r\binom{j+m+r-2}{r}}
\end{equation*}
Now $h(n,1,m)=\sum_jH_1(n,1,m,j)$,
and with $H_1(n,1,m,j)$ as input,
the function {\tt sumrecursion} gives the recurrence
$(n-1)(n-2+m+r)h(n,1,m)-(n+r-1)(n-m)h(n-1,1,m)=0$.
Maple gives the initial value $h(m,1,m)=0$, for
$m>1$.

Case $h(n,n,n)$:
Maple  computes
\begin{equation*}
H_1(n,n,n,j)=\sum_{k=0}^{j-1}H(n,n,n,j,k)=
\frac
{(-1)^{n+j}j\binom{2n+r-2}{n}\binom{n+j+r-3}{j-1}
\binom{n}{j}}
{r\binom{n+j+r-2}{r}}
\end{equation*}
Similarly to the previous case,
{\tt sumrecursion} gives the recurrence
$h(n,n,n)-(n-1,n-1,n-1)=0$,
and obviously $h(1,1,1)=1$.

Case $h(n,i,1)$:
We need to do something different
in this case.
First, we show that our conjectured inverse is
symmetric.
Let 
\begin{multline*}
S(n,i,j,k)=\\
(-1)^{i+k+1}
\binom{n+i+r-2}{i}\binom{n}{i}
\binom{n+k+r-2}{k}\binom{n}{k}
\frac{i^2\prod_{l=0}^{r-3}i+j+l}
{r\prod_{l=0}^{r-2}i+k+l},
\end{multline*}
so that
$B_{ij}(n,r)=\sum_{k=0}^{j-1} S(n,i,j,k)$.
Now {\tt zeil} produces the recurrence 
$S(n+1,i,j,k)-S(n,i,j,k)=T(n,i,j,k+1)-T(n,i,j,k)$,
where
$$
T(n,i,j,k)=\frac{-(2n+r)k^2(i+k+r-2)}{(n+r-1)^2(n-i+1)(n-k+1)}S(n,i,j,k).
$$
This implies that
$B_{ij}(n+1,r)-B_{ij}(n,r)=T(n,i,j,j)-T(n,i,j,0)$.
Now Maple tells us that
$T(n,i,j,j)-T(n,i,j,0)-T(n,j,i,i)+T(n,j,i,0)=0$.
This means that
$B_{ij}(n+1,r)-B_{ji}(n+1,r)=B_{ij}(n,r)-B_{ji}(n,r)$.
Maple tells us that
\begin{multline*}
B_{in}(n,r)-B_{ni}(n,r)=\\
\frac
{\binom{n+i+r-2}{i}i(n+i+r-3)!
\Gamma(2-r)\Gamma(2-n-i-r)}
{r(n+i-1)!(i+r-2)!\Gamma(2-n-r)\Gamma(2-i-r)}
\frac{(-1)^i}{\Gamma(1-i)},
\end{multline*}
which implies $B_{in}(n,r)-B_{ni}(n,r)=0$.

Since $R_n(b_k)$ and $B(n,r)$
are symmetric, the $(1,i)$ entry of
$R_n(b_k)B(n,r)$
equals the $(i,1)$ entry of
$B(n,r)R_n(b_k)$.
The former is
$\sum_{j=1}^n\sum_{k=1}^{i-1}U(n,i,j,k)$,
where
$$
U(n,i,j,k)=
\binom{j+r-1}{r}^{-1}
S(n,j,i,k)
.$$
The function {\tt zeil} produces
$$
Y(n,i,j,k)=\frac{-(2n+r)k^2(j+k+r-2)}{(n+r-1)^2(n-j+1)(n-k+1)}U(n,i,j,k)
$$
which satisfies
$$
U(n+1,i,j,k)-U(n,i,j,k)=Y(n,i,j,k+1)-Y(n,i,j,k).
$$
Thus we have
$$\sum_{k=1}^{i-1}U(n+1,i,j,k)-\sum_{k=1}^{i-1}U(n,i,j,k)=Y(n,i,j,i)-Y(n,i,j,0),$$
and Maple tells us that
$\sum_{j=1}^n Y(n,i,j,i)-Y(n,i,j,0)=0$.
All that remains is to check the initial value
$\sum_{j=1}^i\sum_{k=1}^{i-1}U(i,i,j,k)=0$.
Maple tells us that
\begin{equation*}
\sum_{j=1}^i\sum_{k=1}^{i-1}U(i,i,j,k)=
\frac
{\Gamma(1-r)\Gamma(2i-r)\Gamma(2i+r+1)\Gamma(2i+r-1)}
{\Gamma(-r-1)^2\Gamma(i+r+1)^2\Gamma(i+r)}
\frac
{(-1)^i}
{(i-1)\Gamma(-i)},
\end{equation*}
which implies that
$\sum_{j=1}^i\sum_{k=1}^{i-1}U(i,i,j,k)=0$
when $i>1$.
\end{proof}

We consider reciprocal Hankel matrices
based on one more sequence of binomial
coefficients.
Let $c_k=\binom{k+3}{3}$.

\begin{thm}
\label{b3thm}
Let $C(n)$ be the $n\times n$ matrix whose
$(i,j)$-entry is 
$$
C_{ij}(n)=
\sum_{k=0}^{j-1}
(-1)^{i+k+1}
\binom{n+i+2}{i+k+1}
\binom{n+k+1}{i+k+1}
\binom{i+k+1}{i}
\binom{i+k}{i}
\frac{i(j-k)}{3}.
$$
Then $C_{ij}(n)$ is an integer,
and $C(n)$ is
the inverse of the matrix $R_n(c_k)$.
\end{thm}

\begin{proof}
First we show that each summand of
the sum which defines each entry
is an integer.
It is well known that if $a\equiv 0 \pmod 3$,
$b\equiv 1 \pmod 3$, 
$c\equiv 2 \pmod 3$,
then 
$\binom{a}{b}$, 
$\binom{a}{c}$, and
$\binom{b}{c}3$
are all divisible by 3.
Using this fact, we find that one of the terms
$\binom{i+k+1}{i}$,
$\binom{i+k}{i}$,
or $i$
is divisible by 3 unless
$i\equiv 1 \pmod 3$ and 
$k\equiv 0 \pmod 3$.
But now 
$n+i+2\equiv n \pmod 3$,
$n+k+1\equiv n+1 \pmod 3$, and 
$i+k+1\equiv 2 \pmod 3$.
Thus 3 divides one of the terms
$\binom{n+i+2}{i+k+1}$ or $\binom{n+k+1}{i+k+1}$.

The proof that $C(n)$ is the
inverse of $R_n(c_k)$ is similar to the proof
of Theorem \ref{Mother}.
Let 
$Z(n,i,m,j,k)=
(-1)^{i+k+1}
\binom{n+i+2}{i+k+1}
\binom{n+k+1}{i+k+1}
\binom{i+k+1}{i}
\binom{i+k}{i}
\frac{i(j-k)}{3\binom{j+m+2}{3}}
$,
so that
$z(n,i,m)=\sum_{j=1}^{n}\sum_{k=0}^{j-1}Z(n,i,m,j,k)$
is the $(i,m)$ entry of 
$C(n)R_n(c_k)$.
Then $Z$ satisfies the recurrence
\begin{multline*}
(n-i+1)(n+i+1)\Bigl(Z(n-1,i-1,m,j,k)-Z(n,i-1,m,j,k)\Bigr)+
\\
i(i-1)\Bigl(Z(n-1,i,m,j,k)-Z(n,i,m,j,k)\Bigr)=0.
\end{multline*}
Now the proof proceeds similarly to the 
proof of Theorem \ref{Mother},
except that we don't have to do the
difficult initial value $m=1$.
\end{proof}

One might wonder whether there isn't a simpler
formula
than the one we give for $B_{ij}(n,r)$.
If we fix $i$ and $j$ and consider $B_{ij}$ as a
polynomial
of $n$, then it usually has an irreducible factor
of degree
$\min\{2i-2,2j-2\}$.
Thus it seems unlikely that one could avoid the sum
in the
given formula.
The next section
suggests that the 
given sum is the `right' way to
describe $B_{ij}(n,r)$.

\section{Reciprocal Hankel matrices based on
Fibonomial coefficients}

Remarkably, by changing the exponent of $-1$
and changing the binomial coeffiecients to
Fibonomial coefficients in the formula for
$B_{ij}$,
we get a formula for the entries of the inverses
of reciprocal Hankel matrices based on Fibonomial
coefficients.

Let $d_k=d_k(r)$ be the Fibonomial coefficient
$\Bigl(\binom{k+r-1}{r}\Bigr)$.

\begin{conj}
Let $D(n,r)$ be the $n\times n$ matrix whose
$(i,j)$-entry is
\begin{multline*}
D_{ij}=D_{ij}(n,r)=\sum_{k=0}^{j-1}
(-1)^{e(n,i,k)}
\Biggl(\binom{n+i+r-2}{i}\Biggr)
\Biggl(\binom{n}{i}\Biggr) \times \\
\Biggl(\binom{n+k+r-2}{k}\Biggr)
\Biggl(\binom{n}{k}\Biggr)
\frac{F_i^2\prod_{l=0}^{r-3}F_{i+j+l}}
{F_r\prod_{l=0}^{r-2}F_{i+k+l}}.
\end{multline*}
Then the $D(n,r)$ is
the inverse of the matrix $R_n(d_k)$.
\end{conj}

We have verified this conjecture for
$n\le 16$ and $r\le 10$. (We assume that $r$ is a
positive integer.)
We also observe that
the inverse of a reciprocal Hankel matrix 
based on Fibonomial coefficients
has integer entries exactly when 
the corresponding reciprocal Hankel matrix 
based on binomial coefficients
has integer entries. This may just be
a consequence of known divisibility
properties of the Fibonomials.
It seems likely that this conjecture may
be proved by combining the methods of the
proofs of Theorem \ref{FilbertX} and Theorem
\ref{Mother},
and that it may be extended to 
the corresponding sequence of $x$-Fibonomial
coefficients.

AMS Classification Numbers: 11B39, 11B65, 15A09.
\end{document}